\documentclass[fleqn]{mat01}
\usepackage{times,mathtimy,amssymb,latexsym}
\begin{document}

\setcounter{page}{333} \firstpage{333}

\newcommand{\w}{\wedge}
\newcommand{\f}{\frac}
\newcommand{\s}{\sigma}
\newcommand{\la}{\lambda}
\newcommand{\Si}{\Sigma}
\newcommand{\iy}{\infty}
\newcommand{\what}{\widehat}
\newcommand{\lgra}{\longrightarrow}
\newcommand{\newsection}{\setcounter{equation}{0}}
\renewcommand\theequation{\thesection\arabic{equation}}

\newtheorem{theore}{Theorem}
\renewcommand\thetheore{\arabic{section}.\arabic{theore}}
\newtheorem{theor}[theore]{\bf Theorem}
\newtheorem{rem}[theore]{\it Remark}
\newtheorem{lem}[theore]{\it Lemma}
\newtheorem{propo}[theore]{\rm PROPOSITION}
\newtheorem{definit}[theore]{\rm DEFINITION}

\newcommand{\om}{\omega}%
\newcommand{\bw}{\bigwedge}%
\newcommand{\ehbig}{\emph {\bigwedge}}%
\newcommand{\ehK}{\emph K}%
\newcommand{\ehC}{\emph C}%
\newcommand{\ehGKK}{\emph G//K}%
\newcommand{\ehp}{\emph p}%
\newcommand{\ehGK}{\emph G/K}%
\newcommand{\ehG}{\emph G}%
\newcommand{\R}{\mathbb R}%
\newcommand{\C}{\mathbb C}%
\newcommand{\D}{\mathbb D}%
\newcommand{\Z}{\mathbb Z}%
\newcommand{\Q}{\mathbb Q}%

\def\d{\mbox{\rm d}}
\def\e{\mbox{\rm e}}


\title{On the Schwartz space isomorphism theorem for
rank one symmetric space}

\markboth{Joydip Jana and Rudra P~Sarkar}{Riemannian symmetric space of rank one}

\author{JOYDIP JANA and  RUDRA P SARKAR}

\address{Indian Statistical Institute, Division of Theoretical Statistics and
Mathematics,\break
203 B.T. Road, Kolkata~700~108, India\\
\noindent E-mail: joydip$_{-}$r@isical.ac.in; rudra@isical.ac.in}

\volume{117}

\mon{August}

\parts{3}

\pubyear{2007}

\Date{MS received 9 August 2005; revised 14 May 2007}

\begin{abstract}
In this paper we give a simpler proof of the $L^p$-Schwartz space isomorphism  $(0<
p\leq 2)$ under the Fourier transform for the class of functions of left $\delta$-type
on a Riemannian symmetric space of rank one. Our treatment rests on Anker's \cite{A}
proof of the corresponding result in the case of left $K$-invariant functions on $X$.
Thus we give a proof which relies only on the Paley--Wiener theorem.
\end{abstract}


\maketitle

\section{Introduction}

Let $X$ be a rank one Riemannian symmetric space of noncompact type. We recall that
such a space can be realized as $G/K$, where $G$ is a connected noncompact semisimple
Lie group of real rank one with finite center and $K$ is a maximal compact subgroup of
$G$. Anker \cite{A}, in his paper gave a remarkably short and elegant proof of the
$L^p$-Schwartz space isomorphism theorem for $K$ bi-invariant functions on $G$ under
the spherical Fourier transform for $(0<p \leq 2)$. The result for $K$ bi-invariant
functions was first proved by Harish-Chandra \cite{Ha1,Ha2,Ha3} (for $p=2$) and Trombi
and Varadarajan \cite{TV} (for $0<p<2$). Eguchi and Kowata \cite{M1} addressed the
isomorphism problem for the $L^p$-Schwartz spaces on $X$. In \cite{A}, Anker has
successfully avoided the involved asymptotic expansion of the elementary spherical
functions, which has a crucial role in all the earlier works. In this paper, we have
exploited Anker's technique to obtain the isomorphism of the $L^p$-Schwartz space
$(0<p \leq 2)$ under Fourier transform for functions on $X$ of a fixed $K$-type.

Let $(\delta, V_{\delta})$ be an unitary irreducible
representation of $K$ of dimension $\delta$. Our basic
$L^p$-Schwartz space $S^p_\delta(X)$ is a space of
$\hbox{Hom}(V_\delta, V_\delta)$-valued $C^\infty$ functions, the
Eisenstein integral $\Phi_{\lambda, \delta}(x)$ is a
$\hbox{Hom}(V_\delta, V_\delta)$-valued entire function on $\C$
and $S_\delta(\mathfrak a^*_\C)$ consists of analytic functions on
the strip $\mathfrak a^*_\epsilon =\{\lambda \in \mathbb C
||\hbox{Im}\,\lambda|\leq \epsilon \}$. Anticipating these and
other notations and definitions developed in \S\S2 and 3, we state
the main result of the paper.

\begin{theor}[\!]\label{main-theorem-1}
For $0<~p~\leq 2$ and $\epsilon = 2/p-1$ the $\delta$-spherical
transform $f \mapsto \tilde{f},$\break where
\begin{equation}\label{delta spherical transform}
\tilde{f}(\lambda) = d(\delta) \int_X {\rm tr}\,f(x) \
\Phi_{\bar{\lambda}, \delta} (x)^* {\rm d}x,
\end{equation}
is a topological vector space isomorphism between the spaces $S^p_\delta(X)$ and
$S_\delta(\mathfrak a^*_\epsilon);$ with the inverse
\begin{equation}\label{Inversion of the delta spherical transform}
f(x) = \omega^{-1} \int_\mathfrak {a^*} \ \Phi_{\lambda,
\delta}(x) \tilde{f}(\lambda) |c(\lambda)|^{-2} \ {\rm d}\lambda.
\end{equation}
\end{theor}

\section{Preliminaries}
\setcounter{equation}{0}

The pair $(G, K)$ and $X$ are as described in the introduction. We let $G = KAN$
denote a fixed Iwasawa decomposition of $G$. Let $\mathfrak g$, $\mathfrak k$,
$\mathfrak a$ and $\mathfrak n$ denote the Lie algebras of $G, K, A$ and $N$
respectively. We recall that dimension of $\mathfrak a=1$.

Let $\mathfrak a^*$ be the real dual of $\mathfrak a$ and $\mathfrak a^*_\C$ be its
complexification. We identify $\mathfrak a, \mathfrak a^*$ with $\R$ and $\mathfrak
a^*_\C$ with $\C$ using a normalization explained below. Let $H\hbox{\rm :}\ g \mapsto
H(g)$ and $A\hbox{\rm :}\ g \mapsto A(g)$ be projections of $g \in G$ in $\mathfrak a$
in Iwasawa $KAN$ and $NAK$ decompositions respectively, that is any $g \in G$ can be
written as $g = k \exp H(g) n= n_1 \exp A(g) k_1$. These two are related by $A(g)
=-H(g^{-1})$ for all $g \in G$. Let $M'$ and $M$ respectively be the normalizer and
centralizer of $A$ in $K$. $M$ also normalizes $N$. Let $W = M'/M$ be the Weyl group
of $G$. Here $W= \{ \pm1 \}$. Let us choose and fix a system of positive restricted
roots which we denote by $\Sigma^+$. The real number $\rho$ corresponds to
$\frac{1}{2} \sum_{\alpha \in \Sigma^+} m_\alpha \alpha$ where $m_\alpha$ is the
multiplicity of the root $\alpha$. With a suitable normalization of the basis of $A$
we can identify $\rho$ with $1$. The positive Weyl chamber $\mathfrak a^+ \subset
\mathfrak a$ (${\mathfrak a^*}^+ \subset \mathfrak a^* $) is identified with the
positive real numbers. We denote $x^+$ to be the $\overline{\mathfrak a^+}$ component
of $x \in G$ for the Cartan decomposition  $ G = K\overline{A^+}K=K(\exp~
\overline{\mathfrak a^+}) K$ and let $|x| = x^+$. We have a basic estimate
(Proposition~4.6.11 of \cite{G2}): there is a constant $c > 0$ such that
\begin{equation}\label{Iwasawa<Cartan}
|H(x)| \leq c |x| \ \ \mbox{ for \ $x \in G$ }.
\end{equation}
We note that, any function $f$ on $X$ can also be considered as a function on the group $G$ with the property $f(gk)=f(k)$, where $g \in G$ and $k \in K$.
 Let $x = k a_t k'$ where $a_t = \exp t \in A$, $t \in \mathfrak a \cong  \R$ .
The Haar measure of $G$ for the Cartan decomposition is given by
\begin{equation}\label{Haar measure for Cartan decomposition}
\int_G f(x)\hbox{d}x = \hbox{const} \int_K \hbox{d}k \int_{\mathfrak a^+} \Delta(t)~
\hbox{d}t \int_K \hbox{d}k' f(k a_t k'),
\end{equation}
where $\Delta(t)= \prod_{\alpha \in \Sigma^+} \sinh^{m_\alpha} \alpha(t)$. In the
Iwasawa decomposition, $x = k a_t n$, the Haar measure is
\begin{equation}\label{Haar measure for Iwasawa decomposition}
\int_G f(x)\hbox{d}x = \hbox{const} \int_K \hbox{d}k \int_{\mathfrak a^+}
\hbox{e}^{2t} \hbox{d}t \int_N \hbox{d}n f(k a_t n).
\end{equation}
In both (\ref{Haar measure for Cartan decomposition}) and (\ref{Haar measure for
Iwasawa decomposition}) `const' stands for positive normalizing constants for the
respective cases.

Let $(\delta, V_{\delta})$ be an unitary irreducible
representation of $K$. Let $d(\delta)$ and $\chi_\delta$ stand for
the dimension and character of the representation $\delta$. Let
$V_\delta^M$ be the subspace of $V_\delta$ fixed under
$\delta|_M$; i.e $V_\delta^M = \{ v \in
V_\delta\!\!\mid\!\!\delta(m)v = v, \forall m \in M\}$. Recall
that as $G$ is of real rank one, the dimension of $V_\delta^M $ is
0 or 1 (see \cite{K}). Let $\hat{K}_M$ be the set of all
equivalence classes of irreducible unitary representation $\delta$
of $K$ for which $V_\delta^M \neq \{0\}$. For our result we choose
$\delta \in \hat{K}_M$. We shall also fix an orthonormal basis
$\{v_1, v_2, \dots, v_{d(\delta)}\}$ of $V_{\delta}$ such that
$v_1$ spans $V_\delta^M$.

We shall denote $\mathcal D(X)$ for the space of all $\C$-valued $C^\infty$ functions
on $X$ with compact support. For any function $f \in \mathcal D(X)$, the Helgason
Fourier transform (HFT) (III, \S1 of \cite{He2}) $\mathcal F f$ is defined by
\begin{equation}\label{Helgason Fourier transform}
\mathcal F f(\lambda, kM) = \int_X f(x) \hbox{e}^{(i\lambda -1)H (x^{-1}k)} \
\hbox{d}x.
\end{equation}
Let us fix the notation $\mathcal F f(\lambda, kM) = \mathcal F f(\lambda, k)$. The
inversion formula for HFT for $f \in \mathcal D(X)$ is given by
\begin{equation} \label{HFT Inversion}
f(x) = \frac{1}{\omega} \int_{\mathfrak a^*} \int_{K} \mathcal F f(\lambda, k)
\hbox{e}^{-(i\lambda +1)H(x^{-1}k)} |\textbf{c}(\lambda)|^{-2} \hbox{d} \lambda
\hbox{d}k.
\end{equation}
Here, $\omega = |W|$ is the cardinality of the Weyl group and $\textbf{c}(\lambda)$ is
the Harish-Chandra $\textbf{c}$-function.  For our purpose we shall  need the
following simple estimate on $\textbf{c}(\lambda)$: there exist constants $c,b
> 0$ such that
\begin{equation}
|\textbf{c}(\lambda)|^{-2} ~\leq~c(|\lambda| + 1)^b ~~\mbox{
for $\lambda \in \mathfrak a^*$}
\end{equation}
(see, [IV, Proposition~7.2 of \cite{He3}).

Let $\mathcal D(X , \hbox{Hom}(V_\delta , V_\delta))$ be the space of all $C^\infty$
functions on $X$ taking values in $\hbox{Hom}(V_\delta , V_\delta)$ and with compact
support.

Let $\mathcal D^\delta(X) = \{f \in \mathcal D(X ,
\hbox{Hom}(V_\delta , V_\delta))\!\mid\!f(k \cdot x) = \delta (k)
f(x)\}$.  We topologize $\mathcal D^\delta (X)$ by the inductive
limit topology of the spaces $\mathcal D_R (X ,
\hbox{Hom}(V_\delta , V_\delta))$, where $R = 0,1,2, \dots$. These
are the spaces of functions on $X$ with support lying in the
geodesic $R$-balls. Let $\check{\delta}$ be the contragradient
representation of $\delta$. The class of functions $\mathcal
D_{\check{\delta}} (X) = \{ f \in \mathcal D(X)\!\mid\!f \equiv
d(\delta) \chi_{\delta} \ast f\}$ is the space of all left
$\check{\delta}$ type functions on $X$. Being a subspace of
$\mathcal D (X)$, $\mathcal D_{\check{\delta }}(X)$ inherits the
subspace topology of $\mathcal D (X)$. We also notice that, for $f
\in C^\infty (X)$ the function
\begin{equation}
\label{proj} f^\delta (x) = d(\delta) \int _K f(k\cdot x) \delta (k^{-1}) \hbox{d}k
\end{equation}
is a $C^\infty$ map from $X$ to $\hbox{Hom} (V_\delta , V_\delta)$ satisfying
\begin{equation*}
f^\delta (k \cdot x) = \delta(k) f^\delta (x).
\end{equation*}
The following lemma (III,  Proposition~5.10 of \cite{He2}) shows that the two function
spaces $\mathcal D^\delta(X)$ and $\mathcal D_{\check{\delta}}(X)$ are topologically
isomorphic.

\setcounter{theore}{0}
\begin{lem}\label{Q-map}\hskip -.4pc {\rm [9].} \ \
The map $Q\hbox{\rm :}\ \mathcal D^\delta(X) \longrightarrow \mathcal
D_{\check{\delta}}(X)$ given by
\begin{equation*}
Q\hbox{\rm :}\ f \longmapsto \hbox{\rm tr}\,f
\end{equation*}
is a homeomorphism with the inverse given by $Q^{-1} (g)= g^\delta$ for $g \in
\mathcal D_{\check{\delta}}(X)$.
\end{lem}

\section{The \pmb{$\delta$}-spherical transform}
\setcounter{equation}{0}

Most of the material in this section can be retrieved from \cite{He2}. Here we will
restructure the results in a form which is suitable for our purpose. In particular we
will transfer the results from $D_{\check{\delta}}(X)$ to $\mathcal D^\delta(X)$ using
the homomorphism $Q$, defined in Lemma~\ref{Q-map}.

\setcounter{theore}{0}
\begin{definit}$\left.\right.$\vspace{.5pc}

\noindent {\rm For $f \in \mathcal D^\delta(X)$ the
$\delta$-spherical transform $\tilde{f}$ is given by
\begin{equation}
\label{delta-sp-def} \tilde{f}(\lambda) = d(\delta)\int_X\,
\hbox{tr}\,f(x)\Phi_{\bar{\lambda}, \delta}(x)^* \hbox{d}x,\quad
\lambda \in \C
\end{equation}
where, $\Phi_{{\lambda}, \delta}(x)$ is the generalized spherical function (Eisenstein
integral). Precisely,
\begin{equation}
\label{Pi} \Phi_{{\lambda}, \delta}(x) = \int_K \hbox{e}^{-(i \lambda + 1)H(x^{-1}k)}
\delta(k) \hbox{d}k
\end{equation}
and therefore, the adjoint of $\Phi_{{\lambda}, \delta}(x)$ is
\begin{equation}
\label{Pi*} \Phi_{\bar{\lambda}, \delta}(x)^* = \int_K
\hbox{e}^{(i \lambda - 1)H(x^{-1}k)} \delta(k^{-1}) \hbox{d}k.
\end{equation}}
\end{definit}
The following is a list of some basic properties of the generalized spherical
functions.

\begin{enumerate}\label{M-fixedness of Phi}
\renewcommand\labelenumi{\arabic{enumi}.}
\leftskip -.2pc
\item For $k \in K$, $\Phi_{{\lambda}, \delta}(kx) = \delta(k) \Phi_{{\lambda},
\delta}(x)$  and $\Phi_{{\lambda}, \delta}(kx)^* =
\Phi_{{\lambda}, \delta}(x)^* \delta(k^{-1})$. For $v \in
V_\delta$ and $ m \in M$, $\delta(m) (\Phi_{\bar{\lambda},
\delta}(x)^* v)= \Phi_{\bar{\lambda}, \delta}(x)^* v$. This shows
that $\Phi_{\bar{\lambda}, \delta}^* $ is a $\hbox{Hom}
(V_\delta,V_\delta^M)$-valued function on $X$.

\item Let $\textbf{L}$ be the Laplace--Beltrami operator of $X$. Then
$\textbf{L}~\Phi_{{\lambda}, \delta} = -(\lambda^2
+1)\Phi_{{\lambda}, \delta}$ (\S1(6) of \cite{He2}).

\item Let $\mathcal U(\mathfrak g_{\C})$ be the the universal enveloping algebra of
$G$. For any $g_1 , g_2 \in \mathcal U(\mathfrak g_{\C})$ there exist constants
$c_\delta = c_\delta(g_1 , g_2, \delta), c_0
> 0, b = b(g_1 , g_2) $ so that (see \cite{Ath})
\begin{equation}\label{estimate of gen sp funct}
\hskip -1.25pc \|\Phi_{{\lambda}, \delta}(g_1 , x , g_2)\| \leq
c_\delta (1 + |\lambda|)^b \varphi_0(x) \hbox{e}^{c_0 |{\rm
Im}\,\lambda| (1 + |x|)},~~ x\in X.
\end{equation}
Here $\|\cdot\|$ is the Hilbert--Schmidt norm.

\item If $\delta$ is the trivial representation of $K$ then $\Phi_{{\lambda}, \delta}(x)$
reduces to the elementary spherical function
\begin{equation}\label{elementary spherical function}
\hskip -1.25pc \varphi_\lambda(x) = \int_K \hbox{e}^{-(i \lambda + 1)H(x^{-1}k)}
\hbox{d}k.
\end{equation}
It satisfies the following estimates:

\leftskip -.2pc \vspace{1pc}
(i) For each $H \in \overline{\mathfrak a^+}$ and $\lambda \in
\overline{{\mathfrak a^*}^+}$,
\begin{equation}
\label{estimate of phi lambda} \hskip -2.5pc 0 < \varphi_{-i \lambda}(\exp H) \leq
\hbox{e}^{\lambda H} \varphi_0(\exp H),
\end{equation}
where, $\varphi_0(\cdot)$ is the elementary spherical function at $\lambda = 0$ (see
Proposition~4.6.1 of \cite{G2}).

(ii) For all $g \in G$, $0 < \varphi_0(g) \leq 1$
(Proposition~4.6.3 of \cite{G2}) and for $t \in
\overline{\mathfrak a^+}$,
\begin{equation}\label{estimate of phi_0}
\hskip -2.5pc \hbox{e}^{-t} \leq \varphi_0(\exp t) \leq q (1+t)\hbox{e}^{-t}
\end{equation}
for some $q> 0$ (see \cite{A2} for a sharper
estimate).\vspace{1pc}

\item We have already noticed that $V_\delta^M$ is 1-dimensional. For $\lambda \in \mathfrak a^*_\C$,
$\delta \in \hat{K}_M$ and $x \in X$, the linear functional
$\Phi_{\lambda, \delta}(x)|_{V_\delta^M}$ is a scalar
multiplication. The elementary spherical function $\phi_{\lambda}$
is related to $\Phi_{\lambda, \delta}$ in the following way (see
III, Corollary~5.17 of \cite{He2}):
\begin{equation}\label{Interrelation of spherical and delta spherical}
\hskip -1.25pc \Phi_{\lambda, \delta}(x)\!\mid_{V_\delta^M}= Q^\delta(\lambda)^{-1}
(\textbf{D}^\delta \varphi_{\lambda})(x),
\end{equation}
where $\textbf{D}^\delta$ is a certain constant coefficient differential operator and
$Q^\delta(\lambda)$ is a constant real coefficient polynomial  in $i \lambda$. An
explicit expression for the polynomial $ Q^\delta$ is available in III, \S2 of
\cite{He2}.

\item  For each $a \in A$, the functions $\lambda \mapsto Q^\delta(\lambda)
\Phi_{\lambda, \delta}(a)$ and $\lambda \mapsto
Q^\delta(\lambda)^{-1} \Phi_{\bar{\lambda}, \delta}(a)^*$ are even
holomorphic functions on $\mathfrak a^*_\C$ (see III, Theorem~5.15
of \cite{He2}).
\end{enumerate}

It follows from 1 and 6 above that for $f \in \mathcal
D^\delta(X)$, $\lambda \mapsto
Q^\delta(\lambda)^{-1}\tilde{f}(\lambda)$ is a
$\hbox{Hom}(V_{\delta}$, $V_{\delta}^{M})$-valued even entire
function on $\C$.

The HFT and the $\delta$-spherical transform of a function $f \in \mathcal
D^\delta(X)$ are related in the following manner.

\begin{definit}$\left.\right.$\vspace{.5pc}

\noindent {\rm Let $\delta \in \hat{K}_M$, $f \in \mathcal D(X)$
and $\mathcal F f$ be its HFT. Then let us define the
$\delta$-projection $(\mathcal F f)^\delta$ of $\mathcal F f$ by
\begin{equation}\label{delta-proj of helgason f.t}
(\mathcal F f)^\delta (\lambda , k) = d(\delta) \int_K \mathcal F f (\lambda , k_1k)
\delta(k_1^{-1}) \hbox{d}k_1.
\end{equation}}
\end{definit}

As noted earlier for $f \in \mathcal D(X)$, its $\delta$-projection $f^\delta \in
\mathcal D^\delta(X)$. Each of its matrix entry is a member of $\mathcal D(X)$. We
define the HFT of $f^\delta$  by
\begin{equation}
\mathcal F(f^\delta)(\lambda , k) = \int_X f^\delta(x) \hbox{e}^{(i \lambda -
1)H(x^{-1}k)} \hbox{d}x.
\end{equation}
This is nothing but the usual HFT at each matrix entry of $f^\delta$.

\begin{propo}\label{delta-proj of HFT+property} $\left.\right.$\vspace{.5pc}

\noindent For $f \in \mathcal D(X)$ and $\delta \in \hat{K}_M$ the
following are true$:$
\begin{enumerate}
\renewcommand\labelenumi{{\rm (\arabic{enumi})}}
\leftskip .2pc
\item $(\mathcal F f)^\delta(\lambda, k) = \delta(k)(\mathcal F f)^\delta (\lambda,
e)$.\vspace{.2pc}

\item $\mathcal F (f^\delta) (\lambda, k) = (\mathcal F f)^\delta (\lambda, k)$.
\end{enumerate}
\end{propo}

\begin{proof}
It is clear from the definition that $(\mathcal F f)^\delta(\lambda, k) =
\delta(k)(\mathcal F f)^\delta (\lambda, e)$.

The following straightforward calculation using Fubini's theorem proves the second
assertion.
\begin{align*}
\mathcal F (f^\delta) (\lambda, kM) &= \int_X f^\delta (x) \hbox{e}^{(i
\lambda - 1) H(x^{-1}k)} \hbox{d}x ,\\[.4pc]
&= d(\delta)  \int_X \left\{\int_K f(k_1 x) \delta(k_1^{-1}) \hbox{d}k_1\right\}
\hbox{e}^{(i\lambda - 1) H(x^{-1}k)} \hbox{d}x,
\end{align*}
\begin{align*}
&=d(\delta)\int_K \int_X f(y) \hbox{e}^{(i \lambda - 1) H(y^{-1}k_1 k)}
\delta(k_1^{-1}) \hbox{d}y \hbox{d}k_1,\\[.4pc]
&=d(\delta) \int_X \int_K f(y) \hbox{e}^{(i \lambda - 1) H(y^{-1}k_2 )}
\delta(k) \delta(k_2^{-1}) \hbox{d}k_2 \hbox{d}y,\\[.4pc]
&=d(\delta) \int_K \mathcal F f(\lambda , k_2 M) \delta(k)
\delta(k_2^{-1}) \hbox{d}k_2,\\[.4pc]
&= d(\delta)\int_K \mathcal F f(\lambda , k_3 k M)
\delta(k_3^{-1}) \hbox{d}k_3,\\[.4pc]
&= (\mathcal F f)^\delta (\lambda, kM).
\end{align*}\vspace{-1pc}
\end{proof}

The next lemma relates the $\delta$-spherical transform defined in
(\ref{delta-sp-def}) with the HFT.

\begin{lem} \label{HFT to delta spherical transform}
If $f \in \mathcal D^\delta (X)$ and $\delta \in \hat{K}_M,$
then~~ $\mathcal F f(\lambda ,e)= \tilde{f}(\lambda)$.
\end{lem}

\begin{proof}For any $f \in \mathcal D^\delta (X)$, by Lemma \ref{Q-map},
$f(x) = d(\delta)\int_K\,\hbox{tr}\,f(kx) \delta(k^{-1}) \
\hbox{d}k$. From the definition of HFT (\ref{Helgason Fourier
transform}) we get
\begin{align*}
\mathcal F f(\lambda ,e) &= \int_X f(x) \hbox{e}^{(i \lambda -
1)H(x^{-1})} \hbox{d}x,\\[.4pc]
&=\int_X d(\delta) \int_K\,\hbox{tr}\,f(kx) \delta(k^{-1}) \
\hbox{d}k \ \hbox{e}^{(i \lambda - 1)H(x^{-1})} \ \hbox{d}x.
\end{align*}
Substituting $kx=y$ we have
\begin{align*}
Ff (\lambda, e) &=d(\delta)\int_X\,\hbox{tr}\,f(y)\int_K \
\hbox{e}^{(i \lambda - 1)H(y^{-1}k)}
\delta(k^{-1})\hbox{d}k,\\[.4pc]
&=d(\delta)\int_X\,\hbox{tr}\,f(y)
\Phi_{\bar{\lambda},\delta}(y)^*
\hbox{d}y,\\[.4pc]
&=\tilde{f}(\lambda).
\end{align*}\vspace{-1pc}
\end{proof}

\begin{lem} \label{Inversion and Plancherel of delta spherical Lemma}
The inversion formula for the $\delta$-spherical transform $f
\mapsto \tilde{f}$ is given by the following{\rm :} For each $f\in
\mathcal D^\delta(X),$
\begin{equation}
f(x)=\frac{1}{\omega}\int_{\mathfrak a^*} \Phi_{\lambda,\delta}(x)
\tilde{f}(\lambda) |\hbox{{\rm {\bf c}}}(\lambda)|^{-2} \ {\rm
d}\lambda.
\end{equation}
Moreover$,$
\begin{equation}\label{Plancherel Thm}
\int_X \|f(x)\|^2 \ {\rm d}x = \frac{1}{w} \int_{\mathfrak a^*}
\|\tilde{f}(\lambda)\|^2 |\hbox{{\rm {\bf c}}}(\lambda)|^{-2} \
{\rm d} \lambda.
\end{equation}
Here$,$ the norm $\|\cdot\|$ is the Hilbert--Schmidt norm.
\end{lem}

\begin{proof}
We use the inversion formula for the HFT~(\ref{HFT Inversion}),
Proposition~\ref{delta-proj of HFT+property} and Lemma~\ref{HFT to delta spherical
transform} to obtain
\begin{align*}
f(x)&= \frac{1}{\omega} \int_{\mathfrak a^*} \int_{K} \mathcal F f(\lambda,
k)\hbox{e}^{-(i\lambda +1)H(x^{-1}k)} |\textbf{c}(\lambda)|^{-2} \
\hbox{d} \lambda \hbox{d}k \\[.4pc]
&=\frac{1}{\omega} \int_{\mathfrak a^*} \int_{K} \delta(k)\mathcal F f(\lambda,
e)\hbox{e}^{-(i\lambda +1)H(x^{-1}k)}
|\textbf{c}(\lambda)|^{-2} \hbox{d} \lambda \hbox{d}k\\[.4pc]
&=\frac{1}{\omega} \int_{\mathfrak a^*} \left(\int_{K} \ \hbox{e}^{-(i\lambda
+1)H(x^{-1}k)}\delta(k)\hbox{d}k \right)\mathcal F f(\lambda,
e)|\textbf{c}(\lambda)|^{-2} \ \hbox{d} \lambda\\[.4pc]
&=\frac{1}{\omega}\int_{\mathfrak a^*} \Phi_{\lambda,\delta}(x)
\tilde{f}(\lambda) |\textbf{c}(\lambda)|^{-2} \ \hbox{d}\lambda.
\end{align*}
As the HFT (\ref{Helgason Fourier transform}) of a function $f \in \mathcal
D^\delta(X)$ is defined entry-wise, it is clear that the Plancherel formula for
Helgason Fourier transform is as follows:
\begin{equation}
\label{Plancherel formula for HFT} \int_X \|f(x)\|^2 \hbox{d}x =
\frac{1}{w} \int_{\mathfrak a^*} \int_K \|\tilde{f}(\lambda,
k)\|^2~|\textbf{c}(\lambda)|^{-2} \hbox{d}k \hbox{d} \lambda.
\end{equation}
Using the relation $\tilde{f}(\lambda, k) = \delta(k)
\tilde{f}(\lambda)$ together with the Schur's orthogonality
relation, the formula (\ref{Plancherel Thm}) can be deduced from
(\ref{Plancherel formula for HFT}).
\end{proof}

\begin{definit}$\left.\right.$\vspace{.5pc}

\noindent {\rm A $C^\infty$ function $\psi$ on $\mathfrak a^*_\C$, with values in
$\hbox{Hom}(V_\delta, V^M_\delta)$, is said to be of {\it exponential type} $R$ if
there exists a constant $R \geq 0$ such that for each $N \in \mathbb Z^+$,
\begin{equation*}
\sup_{\lambda \in \mathfrak a^*_\C} \hbox{e}^{-R|{\rm Im}\,\lambda|}(1+|\lambda|)^N
\|\psi(\lambda)\| <+\infty.
\end{equation*}}
\end{definit}

We denote the space of $C^\infty$ function from $\mathfrak a^*_\C \longrightarrow
\hbox{Hom}(V_\delta, V^M_\delta)$ of exponential type $R$ by $\mathcal H^R(\mathfrak
a^*_\C)$. Let $\mathcal H(\mathfrak a^*_\C) = \bigcup_{R>0}\mathcal H^R(\mathfrak
a^*_\C) $. We state the following topological Paley--Wiener theorem for the $K$-types.
The proof of this theorem follows from III, Theorem~5.11 of \cite{He2} and
Lemma~\ref{Q-map}.

\begin{theor}[\!]\label{PW theorem}
The $\delta$-spherical transform defined in Definition~$\ref{delta-sp-def}$ is a
homeomorphism between the spaces $\mathcal D^\delta(X)$ and $\mathcal
P_\delta(\mathfrak a^*_\C),$ where
\begin{equation*}
\mathcal P_\delta(\mathfrak a^*_\mathbb C) = \{ F \in \mathcal
H(\mathfrak a^*_\mathbb C)|( Q^\delta)^{-1} \cdot F \mbox{is an
even entire  function} \}.
\end{equation*}
Here $Q^\delta(\lambda)$ is the polynomial  in $i \lambda$ with
real coefficients introduced in~$(\ref{Interrelation of spherical
and delta spherical})$.
\end{theor}

Let $\mathcal P_0(\mathfrak a^*_\C)$ denote the set of all even functions in $\mathcal
H(\mathfrak a^*_\C)$, with the relative topology. Let $h \in \mathcal P_0(\mathfrak
a^*_\C)$. By definition, $h$ is a $\hbox{Hom}(V_\delta, V_\delta^M)$-valued function.
As $V_\delta^M$ is of dimension~1 so we can write $h=(h_1, \dots, h_{d(\delta)})$,
where each of $h_i$ satisfies the following conditions:

\begin{enumerate}
\renewcommand\labelenumi{(\roman{enumi})}
\leftskip .35pc
\item it is of exponential type,

\item it is entire,

\item it is an even function.
\end{enumerate}
Let $\mathcal D(K \backslash X)$ and $\mathcal D(K \backslash X,
\hbox{Hom}(V_\delta,V_\delta^M) )$ denote the left $K$-invariant, compactly supported,
$C^\infty$ functions on $X$ taking values respectively in $\C$ and
$\hbox{Hom}(V_\delta,V_\delta^M)$. The spherical transform of $\phi \in \mathcal D(K
\backslash X)$ is defined by $\phi \mapsto \int_X \phi(x) \varphi_\lambda(x^{-1})
\hbox{d}x$. For the class $\mathcal D(K \backslash X,\hbox{Hom}(V_\delta,V_\delta^M)
)$ we define it entry-wise. From the Paley--Wiener theorem for the spherical transform
\cite{G2}, there exists one $f_i \in \mathcal D(K \backslash X)$ so that
$h_i(\lambda)= \int_G f_i(x) \varphi_\lambda(x^{-1}) \hbox{d}x$. Therefore $\mathcal
P_0(\mathfrak a^*_\C)$ is the image of $\mathcal D(K \backslash
X,\hbox{Hom}(V_\delta,V_\delta^M) )$ under the spherical transform. The following
lemma shows that the Paley--Wiener (PW) spaces $\mathcal P_\delta(\mathfrak a^*_\C)$
and $\mathcal P_0(\mathfrak a^*_\C)$ are homeomorphic.

\begin{lem}\label{identify PW spaces}\hskip -.4pc {\rm (III, {\it Lemma}~5.12 {\it of} \cite{He2}).} \ \
The mapping
\begin{equation}
\psi(\lambda) \mapsto Q^{\delta}(\lambda)\psi(\lambda)
\end{equation}
is a homeomorphism of $\mathcal P_0(\mathfrak a^*_\C)$ onto $\mathcal
P_\delta(\mathfrak a^*_\C)$.
\end{lem}

\begin{lem}\label{identify delta type and invariant}
Any $f \in \mathcal D^\delta(X)$ can be written as $f(x) =
\hbox{\rm {\bf D}}^\delta \phi(x),$ where $\phi \in \mathcal D(
K\!\!\setminus\!X,$ ${\rm Hom}(V_\delta, V_\delta^M) )$ and
$\hbox{\rm {\bf D}}^\delta $ is a certain  constant coefficient
differential operator.
\end{lem}

\begin{proof}
Let $f \in \mathcal D^\delta(X)$. Then $\tilde{f} \in \mathcal
P_\delta(\mathfrak a^*_\C)$. Therefore by Lemma~\ref{identify PW
spaces}, the map $\lambda \mapsto \Phi(\lambda)=
Q^\delta(\lambda)^{-1}\tilde{f}(\lambda) $ is in $\mathcal
P_{0}(\mathfrak a^*_\C )$. By the PW theorem for the spherical
function we get one $\phi \in \mathcal D(K\!\!\setminus\! X ,
\hbox{Hom}(V_\delta, V_\delta^M))$ such that
\begin{equation}\label{spherical inversion}
\phi(x)= \frac{1}{w} \int_{\mathfrak a^*} \varphi_{\lambda}(x) \Phi(\lambda)
|\textbf{c}(\lambda)|^{-2} \ \hbox{d} \lambda,
\end{equation}
where $\varphi_{\lambda}(\cdot)$ is an elementary spherical function. Now applying the
differential operator $\textbf{D}^\delta$ (see \ref{Interrelation of spherical and
delta spherical})) on both sides of (\ref{spherical inversion}) we get
\begin{align}
(\textbf{D}^\delta \phi)(x) &= \frac{1}{w} \int_{\mathfrak a^*} \Phi_{\lambda, \delta}(x) Q^\delta(\lambda) \Phi(\lambda) |\textbf{c}(\lambda)|^{-2} \ \hbox{d} \lambda
\nonumber,\\[.5pc]
&=\frac{1}{w} \int_{\mathfrak a^*} \Phi_{\lambda, \delta}(x)
\tilde{f}(\lambda) |\textbf{c}(\lambda)|^{-2} \ \hbox{d} \lambda
\nonumber,\\[.5pc]
&= f(x).
\end{align}$\Box $\vspace{-2.5pc}
\end{proof}

We shall denote the Hilbert--Schmidt norm of an operator by $\|H\|$.\vspace{.5pc}

\begin{definit}\label{Schwartz-space-on-X}\hskip -.5pc{\rm (The $L^p$-Schwartz space on $X$).}$\left.\right.$\vspace{.5pc}

\noindent {\rm For every $0<p \leq2$, $\textbf{D},\textbf{ E} \in \mathcal U
(\mathfrak g_\C)$ and $q \in \mathbb N \cup \{0\}$ we define a semi-norm on $f \in
C^{\infty} (X, {\rm Hom}(V_\delta, V_\delta))$ by
\begin{equation}\label{schwartz-on-X}
\nu_{\textbf{D} , \textbf{E} , q}(f) = \sup_{x \in G} \|f (\textbf{D}, x,
\textbf{E})\|~ \varphi_0(x)^{-2/p} (1 + |x|)^q.
\end{equation}
Let $S^p(X)$ be the space of all functions in $C^{\infty} (X, {\rm Hom}(V_\delta,
V_\delta))$ such that $\nu_{\textbf{D} , \textbf{E} , q}(f) < \infty$ for all
$\textbf{D}, \textbf{E} \in \mathcal U (\mathfrak g_\C)$ and $q \in \mathbb N \cup
\{0\}$. We topologize $S^p(X)$ by means of the seminorms $\nu_{\textbf{D} , \textbf{E}
, q,} \textbf{D}, \textbf{E} \in \mathcal U (\mathfrak g_\C)$, $q \in \mathbb N \cup
\{0\}$.}
\end{definit}

Then $S^p(X)$ is a Frechet space and $\mathcal D (X, \hbox{Hom}(V_\delta,V_\delta))$
is a dense subspace of $S^p(X)$. Let $S^p_\delta(X)$ be the subspace of $S^p(X)$
consisting of the left $\delta$ type $\hbox{Hom}(V_\delta,V_\delta)$-valued functions
in $S^p(X)$. Then clearly $\mathcal D^\delta(X)$ is a dense subspace in
$S^p_\delta(X)$.

\begin{rem}\label{Extension-of-Q-map}
{\rm Let $S^p(X)_{\check{\delta}}$ be the Schwartz space of scalar-valued
$\check{\delta}$ type functions. Recall that $\mathcal D_{\check{\delta}}(X)$ is dense
in $S^p(X)_{\check{\delta}}$. Therefore the homeomorphism $Q$ defined in
Lemma~\ref{Q-map} between $\mathcal D^\delta(X)$ and $\mathcal D_{\check{\delta}}(X)$
extends to a homeomorphism between the corresponding Schwartz spaces $S^p_\delta(X)$
and $S^p(X)_{\check{\delta}}$.}
\end{rem}

We shall now define the Schwartz space $S_\delta(\mathfrak a^*_\epsilon)$ containing
the Paley--Wiener space $\mathcal P_\delta(\mathfrak a^*_\mathbb C)$ as follows.

\begin{definit}\label{Schwartz-space-on-lamda}$\left.\right.$\vspace{.5pc}

\noindent {\rm Let $S_\delta(\mathfrak a^*_\epsilon) $ be the
class of functions on $\mathfrak a^*_{\epsilon}$ taking values in
$\hbox{Hom}(V_\delta, V_\delta^M)$ and satisfying the following
conditions:
\begin{enumerate}
\renewcommand\labelenumi{(\arabic{enumi})}
\leftskip .12pc

\item $h$ is analytic in the interior of the strip $\mathfrak a^*_\epsilon$.

\item $h$ extends continuously to the boundary of the strip
$\mathfrak a^*_\epsilon$.

\item $(Q^{\delta})^{-1}~h$ is even and analytic in the
interior of the strip $\mathfrak a^*_\epsilon$.

\item  For each positive integer $r$ and for each symmetric  polynomial
$P$ on $\mathfrak a^*$,
\begin{equation}\label{seminorms on the image of del FT}
\hskip -1.25pc \tau_{r,P}(h) = \sup_{\lambda\,{\rm Int}\,\in \mathfrak a^*_\epsilon }
\|P(\partial \lambda) h(\lambda)\| (1 + |\lambda|)^r <+\infty.
\end{equation}
$P(\partial \lambda)$ is the differential operator obtained by replacing the variable
$\lambda$ by ${\d}/{\d \lambda}$.
\end{enumerate}}
\end{definit}

The topology given by the countable family of seminorms $\tau_{r,P}$ makes
$S_\delta(\mathfrak a^*_\epsilon)$ a Frechet space.

The condition (\ref{seminorms on the image of del FT}) can also be written in the form
\begin{equation*}
\tau_{n,t}(h)= \sup_{\lambda \in\,{\rm Int}\,\mathfrak a^*_\epsilon} \left\|
\left(\frac{\d}{\d \lambda}\right)^t \{(1+\lambda^2)^n h(\lambda)\}\right\| < +\infty.
\end{equation*}

Let $S_0(\mathfrak a^*_\epsilon)$ be the class of all even
functions on $\mathfrak a^*_\C$ taking values in
$\hbox{Hom}(V_\delta,V_\delta^M)$ satisfying conditions (1), (2)
and (4) of Definition~\ref{Schwartz-space-on-lamda}. Then
$S_0(\mathfrak a^*_\epsilon)$ becomes a Frechet space with the
seminorms $\tau_{r,P}$. Clearly, $\mathcal P_0(\mathfrak a^*_\C)
\subset S_0(\mathfrak a^*_\epsilon)$.

\begin{lem}\label{Imp-0}
The map
\begin{equation}
h(\lambda) \mapsto Q^{\delta}(\lambda) h(\lambda)
\end{equation}
is a homeomorphism from $S_0(\mathfrak a^*_\epsilon)$ onto $S_\delta(\mathfrak
a^*_\epsilon)$
\end{lem}

\begin{proof}
Let $h \in S_0(\mathfrak a^*_\epsilon)$. Then
\begin{align*}
&\sup_{\lambda \in\,{\rm Int}\,\mathfrak a^*_\epsilon}\left\| \left( \frac{\d}{\d
\lambda}\right)^t Q^{\delta}(\lambda)h(\lambda)
\right\| (1+|\lambda|)^m \nonumber\\[.4pc]
&\quad\,\leq \sum_{t_i=0}^t c_i \sup_{\lambda \in\,{\rm Int}\,\mathfrak a^*_\epsilon}
\left\|\left(\frac{\d}{\d \lambda}\right)^{t_i} Q^{\delta}(\lambda) \cdot
\left(\frac{\d}{\d
\lambda} \right)^{t-t_i} h(\lambda) \right\|(1+|\lambda|)^m \nonumber,\\[.4pc]
&\quad\,\leq \sum c_i^\delta \sup_{\lambda \in\,{\rm Int}\,\mathfrak a^*_\epsilon}
\left\| \left( \frac{\d}{\d \lambda} \right)^{t-t_i} h(\lambda)\right\|
(1+|\lambda|)^{m_i}.
\end{align*}
The constants $c_i^\delta$ and the positive integers $m_i$ are dependent on $\delta$.
On the other hand, if $g \in S_\delta(\mathfrak a^*_\epsilon)$ then $\psi(\lambda)=
{g(\lambda)}/{Q^{\delta}(\lambda)}$ satisfies  the conditions (1) and (2) of
Definition~\ref{Schwartz-space-on-lamda}. As $g \in S_\delta(\mathfrak a^*_\epsilon)$,
by (3), $\psi$ is an even function. We need to establish (4) of
Definition~\ref{Schwartz-space-on-lamda} to conclude $\psi \in S_0(\mathfrak
a^*_\epsilon)$.

Let us choose a  compact subset \textbf{C} of $\mathfrak a^*_\epsilon$ containing all
the zeros of $Q^{\delta}(\lambda)$ in the strip $\mathfrak a^*_\epsilon$ such that
$|Q^{\delta}(\lambda)| \geq \alpha$ for all $\lambda \in \mathfrak a^*_\epsilon
\setminus \textbf{C} $, where $\alpha$ is a positive constant.
\begin{align*}
&\sup_{\lambda \in\,{\rm Int}\,\mathfrak a^*_\epsilon} \left\|\left(\frac{\d}{\d
\lambda}\right)^t
\psi(\lambda)\right\|(1+|\lambda|)^m\\[.5pc]
&\quad\,\leq \sup_{\lambda \in \textbf{C}}\left\|\left(\frac{\d}{\d \lambda}\right)^t
\frac{g(\lambda)}{Q^{\delta} (\lambda)}\right\|(1+|\lambda|)^m + \sup_{\lambda
\in\,{\rm Int}\,\mathfrak a^*_\epsilon \setminus \textbf{C}} \frac{\|\beta(\lambda)
\big(\frac{{\rm d}}{{\rm
d} \lambda}\big)^{t_1} g(\lambda)\|}{|Q^{\delta}(\lambda)|^{t_2}}\\[.5pc]
&\quad\,\leq k_1 + \frac{k_2}{\alpha} \sup_{\lambda \in\,{\rm Int}\,\mathfrak
a^*_\epsilon} \left\|\left(\frac{\d}{\d \lambda}\right)^{t_1} g(\lambda)\right\|
(1+|\lambda|)^{m_1} < +\infty,
\end{align*}
where $\beta(\lambda)$ is a polynomial in $\lambda$. This concludes the proof. \hfill
$\Box$
\end{proof}

It follows from above that any $h$ in $S_\delta(\mathfrak a^*_\epsilon)$ can be
written as $Q^{\delta} (\lambda) g(\lambda)$ where $g \in S_0(\mathfrak a^*_\epsilon)$
and vice-versa.

Let $g =(g_1,\dots,g_{d(\delta)})\in  S_0(\mathfrak a^*_\epsilon)$. Then each
scalar-valued function $g_i$ belongs to the Schwartz space $S(\mathfrak a^*_\epsilon)$
containing the Paley--Wiener space $\mathcal P(\mathfrak a^*_\C) $ under the spherical
Fourier transform.

\begin{propo}\label{PW-is-dense-in Schwartz-space}$\left.\right.$\vspace{.5pc}

\noindent The Paley--Wiener space $\mathcal P_\delta(\mathfrak a^*_\C)$ is a dense
subspace of $S_\delta(\mathfrak a^*_\epsilon)$.
\end{propo}

\begin{proof}
We have seen in Lemma \ref{Imp-0} that any $h
=(h_1,\dots,h_{d(\delta)}) \in \mathcal P_\delta(\mathfrak
a^*_\C)$ can be written as $Q^{\delta}\cdot
(g_1,\dots,g_{d(\delta)})$, where each $g_i$ belongs to the
Paley--Wiener space $\mathcal P(\mathfrak a^*_\C)$ under the
spherical Fourier transform. We recall that $\mathcal P(\mathfrak
a^*_\C)$ is dense in $S(\mathfrak a^*_\epsilon)$. Let
$H=(H_1,\dots,H_{d(\delta)}) \in S_\delta(\mathfrak
a^*_\epsilon)$. Then $H=
(Q^{\delta}G_1,\dots,Q^{\delta}G_{d(\delta)} )$ where $G= (G_1,
\dots,G_{d(\delta)})$ $\in S_0(\mathfrak a^*_\epsilon)$, i.e.,
each $G_i \in S(\mathfrak a^*_\epsilon)$. Then there exists a
sequence $\{{g_i}_n\}$ in $\mathcal P_0(\mathfrak a^*_\C)$
converging to $G_i$. Hence $\{Q^{\delta}\cdot {g_i}_n\}$ converges
to $Q^{\delta}\cdot G_i $ in the topology of $S(\mathfrak
a^*_\epsilon)$. Therefore, the sequence $\{Q^{\delta}\cdot
({g_1}_n,\dots,{g_{d(\delta)}}_n)\} \subset \mathcal P_\delta
(\mathfrak a^*_\C)$ converges to
$(Q^{\delta}G_1,\dots,Q^{\delta}G_{d(\delta)} )$ in
$S_\delta(\mathfrak a^*_\epsilon)$. This completes the
proof.\hfill $\Box$
\end{proof}

\section{Proof of Theorem~\ref{main-theorem-1}}

\setcounter{equation}{0}

\setcounter{theore}{0}
\begin{lem}\label{Imp Lemma-1}
Let $f \in S^p_\delta(X)$. Then its $\delta$-spherical transform
$\tilde{f}$ is an analytic function in the interior of the strip
$\mathfrak a^*_\epsilon$.
\end{lem}

\begin{proof}
For any function $f\hbox{\rm :}\ X \mapsto \hbox{Hom}(V_\delta,
V_\delta^M)$, it is easy to show that $|\hbox{tr}\,f(x)| \leq
\|f(x)\|$ for all $x \in X$. As $f \in S^p_\delta(X)$, from
(\ref{schwartz-on-X}), we conclude that for each $\textbf{D},
\textbf{E} \in \mathcal U(\mathfrak g_\C)$ and $n \in \Z^+ \cup
\{0\}$,
\begin{equation}
\label{decay of Tr f} \sup_{x \in X} \|\hbox{tr}\,f(\textbf{D}, x,
\textbf{E})\| (1 + |x|)^n \varphi_0^{-{2}/{p}}(x) < +\infty.
\end{equation}
Using (\ref{decay of Tr f}) and the estimate (\ref{estimate of gen sp funct}) one can
show that the integral in Definition~\ref{delta-sp-def} of the $\delta$-spherical
transform converges absolutely for $\lambda \in \mathfrak a^*_\epsilon$.

A standard application of Morera's theorem together with Fubini's
theorem shows that $\lambda \mapsto \tilde{f}(\lambda)$ is
analytic in the interior of the strip $\mathfrak a^*_\epsilon$.
\hfill $\Box$
\end{proof}

\begin{lem}\label{Imp Lemma-2}
For $f \in S^p_\delta(X)$ and for each $t,n \in \Z^+ \cup \{0\}${\rm ,} there exists a
positive integer $m$ and $n$ such that
\begin{align*}
\sup_{\lambda \in\,{\rm Int}\,\mathfrak a^*_\epsilon}
\left\|\left(\frac{\d}{\d \lambda}\right)^t \{(1+\lambda^2)^n
\tilde{f}(\lambda)\}\right\| \leq c \sup_{x \in X} \|\hbox{\rm
{\bf L}}^{n}f(x)\| (1+|x|)^{m} \varphi_0^{-{2}/{p}}(x),
\end{align*}
where $c$ is a positive constant.
\end{lem}

\begin{proof}
From (\ref{delta-sp-def}) we have
\begin{align}\label{inj1}
\left(\frac{\d}{\d \lambda}\right)^t \{(1+\lambda^2)^{n}
\tilde{f}(\lambda) \} &= \left(\frac{\d}{\d\lambda} \right)^t
\left\{d(\delta) \int_X\,\hbox{tr}\,f(x) (1
+\lambda^2)^{n}\Phi_{\bar{\lambda}, \delta}(x)^* \d
x\right\}\nonumber\\[.4pc]
&= \left(\frac{\d}{\d\lambda}\right)^t \left\{
d(\delta)\int_X\,\hbox{tr}
f(x)(-\textbf{L})^{n}\Phi_{\bar{\lambda}, \delta}(x)^* \d x
\right\},
\end{align}
where the last equality follows from (2) of the discussion following Definition~3.1.
Using integration by parts we get from above that
\begin{align*}
&\left(\frac{\d}{\d\lambda}\right)^{t} \{(1 +
\lambda^{2})^{n}\tilde{f}(\lambda)\}\\[.4pc]
&\quad\,=\left(\frac{\d}{\d\lambda}\right)^t \left\{
d(\delta)\int_X
(-\textbf{L})^{n}\,\hbox{tr}\,f(x)\Phi_{\bar{\lambda},
\delta}(x)^* \d x\right\}\\[.4pc]
&\quad\,=d(\delta) \int_X (-\textbf{L})^{n}\,\hbox{tr}\,f(x)
\left(\frac{\d}{\d\lambda} \right)^t \int_K \e^{(i \lambda -
1)H(x^{-1}k)} \delta(k^{-1}) \d k \d x\\[.4pc]
&\quad\,=d(\delta) \int_X (-\textbf{L})^{n}\,\hbox{tr}\,f(x)~
\int_K (i
H(x^{-1}k))^t~\e^{(i \lambda - 1)H(x^{-1}k)} \delta(k^{-1}) \d k \d x\\
&\quad\,= (i)^{t'}d(\delta) \int_X \int_K
(H(x^{-1}k))^t~\textbf{L}^{n}\,\hbox{tr}\,f(x)~\e^{(i \lambda -
1)H(x^{-1}k)}
\delta(k^{-1}) \d k \d x \\[.4pc]
&\quad\,= (i)^{t'}d(\delta) \int_K \int_X (H
(x^{-1}k))^t~\textbf{L}^{n}\,\hbox{tr}\,f(x)~\e^{(i \lambda -
1)H(x^{-1}k)} \delta(k^{-1}) \d x \d k.
\end{align*}
We substitute $x^{-1}k = y^{-1}$ and use
$\textbf{L}\,\hbox{tr}\,f(y) = \hbox{tr}\,(\textbf{L}\,f)(y)$ to
obtain
\begin{align*}
&\left(\frac{\d}{\d\lambda}\right)^{t} \{(1 +
\lambda^{2})^{n}\tilde{f}(\lambda)\}\\
&\quad\,= (i)^{t'}
d(\delta)\int_{X}\int_{K}(H(y^{-1}))^{t}\,\hbox{tr}\,(\textbf{L}^{n}
f)(ky)~\e^{(i \lambda - 1)H(y^{-1})} \delta(k^{-1}) \d k \d y.
\end{align*}
Note that $\textbf{L}^{n}f$ is again a function of left $\delta$
type. Therefore from above we get
\begin{align}
&\left(\frac{\d}{\d \lambda}\right)^t \{(1+\lambda^2)^n
\tilde{f}(\lambda) \}\nonumber\\[.4pc]
&\quad\,= (i)^{t'}\int_X H (y^{-1})^t~\e^{(i \lambda -
1)H(y^{-1})}\left\{
d(\delta)\int_K\,\hbox{tr}\,(\textbf{L}^{n}f)(ky)\delta(k^{-1})\d
k \right\} \d y,
\nonumber\\[.4pc]
&\quad\,= (i)^{t'}\int_X (H (y^{-1}))^t \textbf{L}^{n}f(y)~\e^{(i
\lambda - 1)H(y^{-1})} \d y\quad \mbox{(by (\ref{proj}))}\nonumber\\[.4pc]
&\quad\,=(i)^{t'}\int_X (H (y))^t \textbf{L}^{n}f(y^{-1})~\e^{(i \lambda - 1)H(y)} \d
y. \label{inj2}
\end{align}
We use the Iwasawa decomposition $G = KAN$ and write $y = ka_{r}n$, where $r\in
{\mathfrak a}$ and $\exp r = a_{r}$ to obtain
\begin{align}
&\left(\frac{\d}{\d \lambda}\right)^t \{(1+\lambda^2)^n
\tilde{f}(\lambda) \}\nonumber\\[.4pc]
&\quad\,= c(i)^{t'} \int_K \int_{\mathfrak a} \int_N \textbf{L}^{n}f(n^{-1} a_r^{-1}
k^{-1}) (H( k a_r n))^t \e^{(i
\lambda - 1)H(k a_r n)} \d k \e^{2r} \d r \d n\nonumber\\[.45pc]
&\quad\,= (i)^{t'} \int_{\mathfrak a} \int_N  \textbf{L}^{n}f((a_r
n)^{-1})~ r^t \e^{(i \lambda+1)r} \d r\d n.
\label{inj3} 
\end{align}
From (\ref{inj3}) we get the following norm inequality.
\begin{align}
\bigg\| \left(\frac{\d}{\d \lambda}\right)^t& \{(1+\lambda^2)^n
\tilde{f}(\lambda) \} \bigg\| \leq c \int_{\mathfrak a} \int_N
\|\textbf{L}^{n}f((a_r n)^{-1})\|~ |r|^t \e^{(| {\rm
Im}\,\lambda|+1)r} \d r\d n.\label{inj4}
\end{align}
As $f \in S^p_\delta(X)$, for each $m \in \Z^+$ we have
$\|\textbf{L}^n f((a_r n)^{-1})\| \leq \nu_{\textbf{L}^n, m} (f)
(1+ |(a_r n)^{-1}|)^{-m}$ $\varphi_0^{{2}/{p}}((a_r n)^{-1})$
where $\nu$ is as defined in (\ref{schwartz-on-X}). Using
(\ref{Iwasawa<Cartan}) we get from above
\begin{align}
&\hskip -4pc \left\|\left(\frac{\d}{\d\lambda}\right)^{t}\{(1 +
\lambda^{2})^{n}\tilde{f}(\lambda)\}
\right\|\nonumber\\[.4pc]
&\hskip -4pc \quad\leq c \nu_{\textbf{L}^n, m} (f) \int_{\mathfrak a} \int_N (1+ |(a_r
n) |)^{-m}\varphi_0^{{2}/{p}}((a_r n)^{-1})(1+ |r|)^t \e^{(|
{\rm Im}\,\lambda|+1)r} \d r\d n\nonumber\\[.4pc]
&\hskip -4pc \quad\leq c_{1} \nu_{\textbf{L}^n, m} (f) \int_{\mathfrak a} \int_N (1+
|(a_r n) |)^{-m+t}\varphi_0^{{2}/{p}}((a_r n)^{-1}) \e^{(| {\rm
Im}\,\lambda|+1)H(a_{r}n)} \d r\d n\nonumber\\[.4pc]
&\hskip -4pc \quad= c_1 \nu_{\textbf{L}^n, m} (f) \int_G (1+ |x|)^{t-m}
\varphi_0^{{2}/{p}}(x) \e^{(|{\rm Im}\,\lambda|-1)H(x)} \d x.\label{inj4}
\end{align}
For convenience we denote $c_1 \nu_{\textbf{L}^n, m} (f)$ by $c_{\nu}$. We use the
Cartan decomposition $G = K\overline{A^{+}}K$ and write $x = k_1 \exp|x| k_2$ and
decompose the integral (\ref{inj4}) as follows:\pagebreak
\begin{align*}
&c_\nu \int_K \int_{\mathfrak a^+} \int_K (1+ |k_1 \exp|x|
k_2|)^{-m+t} \varphi_0^{{2}/{p}}(\exp|x^{-1}|)\\[.5pc]
&\qquad\ \times \e^{(|{\rm Im}\,\lambda|-1)H(\exp|x| k_2)}
\d k_1 \Delta(|x|)\d|x| \d k_2\\[.5pc]
&\quad =c_\nu \int_{\mathfrak a^+} \int_K (1+|x|)^{-m+t}
\varphi_0^{{2}/{p}}(\exp|x|)\\[.5pc]
&\qquad\ \times \e^{(|{\rm Im}\,\lambda|-1)H(\exp|x| k_2)} \Delta(|x|)\d|x| \d k_2,
\end{align*}
as $|x^{-1}|=|x|$ and $|k_1 \exp|x| k_2|=|x|$. Using (\ref{estimate of phi lambda})
the expression above is
\begin{align}
&=c_\nu \int_{\mathfrak a^+} (1+|x|)^{-m+t} \varphi_0^{{2}/{p}}(\exp|x|)\nonumber\\
&\qquad\ \times \left\{\int_K \e^{(-i(i|{\rm Im}\,\lambda|)-1)H(\exp|x| k_2)} \d k_2
\right\}
\Delta(|x|)\d|x|\nonumber\\[.5pc]
& =c_\nu \int_{\mathfrak a^+} (1+|x|)^{-m+t} \varphi_0^{{2}/{p}}(\exp|x|)
\varphi_{-i|{\rm Im}\,\lambda|}(\exp|x|)~\Delta(|x|)\d|x|\nonumber\\[.5pc]
&\leq c_\nu \int_{\mathfrak a^+} (1+|x|)^{-m+t}
\varphi_0^{({2}/{p})+1}(\exp|x|)~\e^{|{\rm
Im}\,\lambda||x|}~\Delta(|x|)\d|x|.\label{inj5}
\end{align}
We take $\lambda\in \hbox{Int}\,\mathfrak a^*_\epsilon$, i.e., $|{\rm Im}\,\lambda| <
\epsilon = \big(\frac{2}{p}-1\big)$. Using the estimate (\ref{estimate of phi_0}) we
get
\begin{align}
&\leq c_\nu \int_{\mathfrak a^+}  (1+|x|)^{(-m+t+ \frac{2}{p}-1)}
\varphi_0^{2}(\exp|x|) ~\Delta(|x|)\d|x|
\nonumber\\[.5pc]
&\label{eq4.8}\leq c_\nu \int_G (1+|x|)^{(-m+t+
\frac{2}{p}-1)}\varphi_0^{2}(x) \d x,\quad (\hbox{see} \
(\ref{Haar measure for Cartan decomposition})).
\end{align}\vspace{-1pc}
\end{proof}
Choosing a suitably large $m$, we see that the integral in
(\ref{eq4.8}) converges (Lemma~11 of \cite{Ha3}). Hence, we
conclude that
\begin{equation}
\sup_{\lambda \in\,{\rm Int}\,\mathfrak a^*_\epsilon}
\left\|\left(\frac{\d}{\d \lambda}\right)^t \{(1+\lambda^2)^n
\tilde{f}(\lambda)\}\right\| \leq \hbox{const} \
\nu_{\textbf{L}^n, m}(f).
\end{equation}
This completes the proof of the lemma. \hfill $\Box$

\begin{lem}\label{conti}
The $\delta$-spherical transform $f \mapsto \tilde{f}$ is a
continuous injection of $S^p_\delta(X)$ into $S_\delta(\mathfrak
a^*_\epsilon)$.
\end{lem}

\begin{proof}
From Lemma~\ref{Imp Lemma-1}, Lemma~\ref{Imp Lemma-2} and (6) of
the discussion below Definition~3.1, we conclude that if $f \in
S^p_\delta(X)$ then $\tilde{f} \in S_\delta(\mathfrak
a^*_\epsilon)$. Also the transform $f\mapsto \tilde{f}$ is
continuous. The fact that $f \mapsto \tilde{f}$ is injective is a
consequence of the Plancherel formula for the HFT (III,
Theorem~1.5 of \cite{He2}).\hfill $\Box$
\end{proof}

The next lemma is an extension of the inversion formula given in Lemma~3.5 for the
Schwartz class functions.

\begin{lem}
Let $h \in S_\delta(\mathfrak a^*_\epsilon)$. Then the inversion $\mathcal Ih$ given
by
\begin{equation*}
\mathcal I h(x) = \frac{1}{\omega} \int_{\mathfrak a^*} \Phi
_{\lambda , \delta}(x) h(\lambda) |\hbox{{\rm {\bf
c}}}(\lambda)|^{-2}\d \lambda
\end{equation*}
is a left $\delta$-type $C^\infty$ function on $X$ taking values in $\hbox{\rm
Hom}(V_\delta, V_\delta)$.
\end{lem}

\begin{proof}
Let us take any derivative $\textbf{D}$ of $X$. For any $\textbf{D} \in \mathcal{U}
({\mathfrak g}_{\mathbb{C}})$,
\begin{equation}\label{*****}
\mathcal Ih(\textbf{D};x)=\frac{1}{\omega} \int_{\mathfrak a^*} \Phi _{\lambda ,
\delta}(\textbf{D};x) h(\lambda) |\textbf{c}(\lambda)|^{-2}\d \lambda.
\end{equation}
Therefore,
\begin{align*}
\|\mathcal Ih(\textbf{D};x)\| &\leq c \int_{\mathfrak a^*} \|\Phi _{\lambda ,
\delta}(\textbf{D};x)\|\|h(\lambda)\|(1+|\lambda|)^b\d\lambda\\[.4pc]
&\leq c_\delta \int_{\mathfrak a^*}
(1+|\lambda|)^{b_{\textbf{D}}+b-n} \varphi_0(x) \d\lambda,\\[.4pc]
&\quad\,\,\begin{array}{@{}l}\hbox{(using estimate (3) of the discussion below}\\
\hbox{Definition~3.1 and (\ref{seminorms on the image of del FT}))}\end{array}\\[.4pc]
& \leq c_\delta \int_{\mathfrak a^*} (1+|\lambda|)^{b_{\textbf{D}}+b-n}\d\lambda.
\end{align*}
We choose $n$ sufficiently large so that the last integral on the right-hand side
exists. Hence, $\mathcal I h (\textbf{D}; x)$ exists for every $\textbf{D}$. Therefore
$\mathcal I h $ is a $C^\infty$ function on $X$. As $\Phi _{\lambda , \delta}(x)$ is
of left $\delta$ type, so is $\mathcal I h$.\hfill $\Box$
\end{proof}

\begin{lem}\label{Surjectivity of delta FT}
If $h \in S_\delta(\mathfrak a^*_\epsilon)${\rm ,} then $\mathcal I h \in S^p_\delta
(X)$.
\end{lem}

\begin{proof}
We consider the spaces $\mathcal P_\delta(\mathfrak a^*_\C)$ and $\mathcal
D^\delta(X)$ equipped with the topologies of the Schwartz spaces $S_\delta(\mathfrak
a^*_\epsilon)$ and $S^p_\delta(X)$ respectively. It is clear from the Paley--Wiener
theorem that $\mathcal I$ maps $\mathcal P_\delta(\mathfrak a^*_\C)$ onto $\mathcal
D^\delta(X)$. We shall show that $\mathcal I$ is a continuous map from $\mathcal
P_\delta(\mathfrak a^*_\C)$ onto $\mathcal D^\delta(X)$ in these topologies. Let $h
\in \mathcal P_\delta(\mathfrak a^*_\C) $ and $\mathcal Ih=f \in \mathcal
D^\delta(X)$. We have to show that for any seminorm $\nu$ on $\mathcal D^\delta(X)$
there exists a seminorm $\tau$ on $P(\mathfrak a^*_\epsilon)$ so that
\begin{equation*}
\nu(f) \leq c_\delta \tau(h),
\end{equation*}
where $c_\delta$ is a constant depending only on $\delta$.

Let $\textbf{D} \in \mathcal U(\mathfrak g_\C)$ and $n \in \Z^+$.
We consider $f$ as a right $K$-invariant function on the group
$G$. Let
\begin{equation}\label{!!!!!!}
\nu_{\textbf{D},n} (f) =\sup_{x\in G} \|\textbf{D} f(x)\|(1+|x|)^n
\varphi_0^{-{2}/{p}}(x).
\end{equation}
From Lemmas~\ref{identify PW spaces} and \ref{identify delta type
and invariant} we know that $f(x)= \textbf{D}^{\delta} \phi(x)$,
where $\phi$ is a $K$ bi-invariant function on $G$ and
$h(\lambda)=Q^{\delta}(\lambda)\Phi(\lambda)$. Here $\Phi$ is the
spherical Fourier transform of $\phi$. Hence from (\ref{!!!!!!})
we have
\begin{equation}
\nu_{\textbf{D},n} (f) =\sup_{x\in G} \|\textbf{D}
\textbf{D}^{\delta} \phi(x)\|(1+|x|)^n
\varphi_0^{-{2}/{p}}(x)=\nu_{\textbf{D}\textbf{D}^{\delta},n}
(\phi).
\end{equation}
By the isomorphism of the $K$ bi-invariant functions in the
Schwartz space (see \cite{A}), for each $\textbf{D} \in \mathcal
U(\mathfrak g), \textbf{D}^{\delta} \in \mathcal U(\mathfrak g)$
and for each $n \in \Z^+$ there exists $m_\delta,t \in \Z^+$ and a
positive constant $c_\delta$ so that,
\begin{equation}
\sup_{x\in G} \|\textbf{D} \textbf{D}^{\delta} \phi(x)\|(1+|x|)^n
\varphi_0^{-{2}/{p}}(x) \leq c_\delta  \sup_{\lambda \in\,{\rm
Int}\,\mathfrak a^*_\epsilon} \left\|\left(
\frac{\d}{\d\lambda}\right)^t \Phi(\lambda)\right\|\!
(1+|\lambda|)^{m_\delta}.
\end{equation}
Now by Lemma~\ref{Imp-0}, for $t,m_\delta \in \Z^+$ there exists
$t_1,m_1 \in \Z^+$ such that
\begin{align*}
\sup_{\lambda \in \mathfrak a^*_\epsilon} \left\|\left( \frac{\d}{\d\lambda}\right)^t
\Phi(\lambda)\right\| (1+|\lambda|)^m &\leq c'_\delta \sup_{\lambda \in\,{\rm
Int}\,\mathfrak a^*_\epsilon} \left\|\left( \frac{\d}{\d\lambda}\right)^{t_1}
h(\lambda)\right\| (1+|\lambda|)^{m_1}
\\[.4pc]
&= c'_\delta \tau_{t_1, m_1}(h) < +\infty .
\end{align*}
Hence, $\nu_{\textbf{D},n} (f) \leq c'_\delta c_{\delta}\tau_{t_1,
m_1}(h)$. The positive constants $c_\delta$ and ${c}'_\delta$ are
dependent on $|\delta|$. The positive integer $m_1$ can be made
independent of the $\delta \in \hat{K}_M$ chosen. This shows that
the inversion $\mathcal I$ is a continuous linear transformation
on a dense subset $\mathcal P_\delta(\mathfrak a^*_\C)$ of
$S_\delta(\mathfrak a^*_\epsilon)$ onto $\mathcal D^\delta(X)$
(The surjectivity follows from Theorem~\ref{PW theorem}.)

Let us now take $h \in S_\delta(\mathfrak a^*_\epsilon) $. As
$\mathcal P_\delta(\mathfrak a^*_\C)$ is dense in
$S_\delta(\mathfrak a^*_\epsilon)$, there exists a Cauchy sequence
$\left\{h_n \right\} \subset \mathcal P_\delta(\mathfrak a^*_\C)$
converging to $h$. Then by what we have proved above, we can get a
Cauchy sequence $\left\{ f_n\right\} \subset \mathcal D^\delta(X)$
such that $\tilde{f_n} = h_n$. As $S^p_\delta(X)$ is a Frechet
space, the sequence converge to some $f \in S^p_\delta(X)$.
Clearly, $f = \mathcal I h $. This completes the\break
proof.\hfill $\Box$
\end{proof}

Finally, Lemmas~\ref{PW theorem}, \ref{conti} and
\ref{Surjectivity of delta FT} together show that the
$\delta$-spherical transform is a surjection onto
$S_\delta(\mathfrak a^*_\epsilon)$ and that $\mathcal I\hbox{\rm
:}\ S_\delta(\mathfrak a^*_\epsilon) \rightarrow S^p_\delta(X) $
is continuous. That is, the $\delta$-spherical transform is a
topological isomorphism between the spaces $S^p_\delta(X)$ and
$S_\delta(\mathfrak a^*_\epsilon)$. \,This proves
Theorem~\ref{main-theorem-1}.

\end{document}